\documentclass[conference]{IEEEtran}
\usepackage{amsmath,amssymb,amsfonts}
\usepackage{algorithm}
\usepackage{algpseudocode}
\usepackage{graphicx}
\usepackage{xcolor}
\usepackage{cite}
\usepackage{booktabs}
\usepackage[most]{tcolorbox}
\usepackage{titlesec}
\usepackage{amsmath}
\usepackage{amssymb}
\usepackage{amsthm} 

\algrenewcommand\algorithmicindent{0.9em}
\algnewcommand\algorithmicparallel{\textbf{in parallel}}

\titlespacing\section{0pt}{4pt plus 4pt minus 2pt}{4pt plus 2pt minus 2pt}
\titlespacing\subsection{0pt}{4pt plus 4pt minus 2pt}{4pt plus 2pt minus 2pt}
\titlespacing\subsubsection{0pt}{4pt plus 4pt minus 2pt}{4pt plus 2pt minus 2pt}
\titlespacing\paragraph{0pt}{4pt plus 4pt minus 2pt}{4pt plus 2pt minus 2pt}
\titleformat*{\subsubsection}{\bfseries}
\titleformat*{\paragraph}{\small\bfseries}

\newtheorem{assumption}{Assumption}
\newtheorem{proposition}{Proposition}
\newtheorem{corollary}{Corollary}
\newtheorem{lemma}{Lemma}
\newtheorem{remark}{Remark}

\IEEEoverridecommandlockouts

\begin{document}

\title{Distributed Time-Varying Resource Allocation in Open Networks}\vspace{-0.2cm}


\author{Thinh T. ~Doan
\thanks{This work is supported in part by NSF under CAREER Award 2527059 and AFOSR Grant FA9550-25-1-0247. Thinh T. Doan is with the University of Texas at Austin, Austin,
TX, 78712 USA (e-mail:{thinhdoan@utexas.edu})}
}
\maketitle

\begin{abstract}
We consider the problem of distributed resource allocation over a team of agents whose active membership changes over time. Each agent’s allocation evolves according to its own temporal dynamics while a shared equality constraint couples all active agents at every instant. For the case of quadratic local objectives, we show that the optimal solution is characterized by a linear system whose coefficient matrix couples a block-tridiagonal temporal structure with a block-diagonal Laplacian capturing the time-varying communication graph. Building on this structure, we develop a distributed forward--backward sweeping algorithm in which each agent resolves its own temporal coupling through a local forward--backward sweep, while the shared resource constraint is resolved through a distributed consensus iteration that is feasible at every iteration. Under a mild diagonal-dominance condition on the temporal coupling, we establish that the optimal allocation at any fixed time converges geometrically to a stationary solution as the horizon grows, and that our distributed algorithm tracks this limit at a matching geometric rate. These results justify a memory-efficient variant of the algorithm in which each agent retains only a bounded window of past history, with the resulting approximation error decaying exponentially in the window size. We validate these theoretical guarantees through numerical simulation, confirming the predicted convergence rates under dynamically changing team membership.
\end{abstract}

\section{Introduction}
We consider distributed resource allocation over an open multi-agent network, where the number and identity of agents may vary over time. Let $t = 0, 1, \ldots, T$ denote a finite time horizon, and let $\mathcal{V}_{t}$ denote the set of agents active at time $t$, with cardinality $\mathcal{N}_t = |\mathcal{V}_{t}|$. Each active agent $i\in\mathcal{V}_{t}$ maintains a decision variable $x_{i,t}\in\mathbb{R}$ and a local objective function $f_{i,t}$. Let $x_t = (x_{i,t})_{i\in\mathcal{V}_t}$. We seek to solve the following time-varying optimization problem
\begin{equation}
\min_{x_1,\dots,x_T} \sum_{t=1}^{T} f_t(x_{t-1},x_t)
\quad \text{s.t.} \quad \sum_{i\in\mathcal{V}_t} x_{i,t} = b_t,\; \forall t,
\label{problem:main}
\end{equation}
where $b_t\in\mathbb{R}$ is a (possibly time-varying) shared resource and $f_{t}$ aggregates the local objectives of active agents at time $t$,
\begin{equation}
f_t(x_{t-1},x_t) = \sum_{i\in\mathcal{V}_t} f_{i,t}(x_{i,t-1},x_{i,t}).\label{problem:obj}
\end{equation}
In our setting, agents may join or leave the network at any time step, so that $\mathcal{V}_{t}$ and $\mathcal{N}_t$ are time-varying, and the dimension of $x_t$ varies accordingly. At each time $t$, agents can communicate with their neighbors characterized by a graph $\mathcal{G}_{t} = (\mathcal{V}_{t},\mathcal{E}_{t})$. Crucially, each $f_{i,t}$ depends on both the current decision $x_{i,t}$ and the previous one $x_{i,t-1}$, coupling every agent’s successive time steps. 

Problem \eqref{problem:main} is difficult for two compounding reasons. First, the temporal coupling in each $f_{i,t}$ prevents decomposing \eqref{problem:main} into independent subproblems, since an agent’s optimal allocation at $t$ depends on its allocation at $t-1$ and on the shared constraint across all active agents. Second, the problem grows with the horizon: increasing $T$ appends a new cost term and variable block, which  alter the optimal solution at every earlier time step. Resolving the entire solution is inefficient as it requires to store the local objectives of all the agents. This raises the central question addressed in this paper: \vspace{-0.2cm}
\begin{tcolorbox}[colback=blue!15!white,colframe=white,lowerbox=invisible]
\textbf{Under what conditions, problem \eqref{problem:main} can be solved efficiently with limited memory?} \vspace{-0.2cm}\end{tcolorbox}\vspace{-0.1cm}


\noindent\textbf{Main Contributions.} In this paper, we study problem \eqref{problem:main} when each agent’s local objective is quadratic, and show that its optimal solution is characterized by a linear system whose coefficient matrix couples a block-tridiagonal temporal structure with a block-diagonal Laplacian capturing the time-varying communication graphs. Building on this structure, we develop a distributed forward-backward sweeping algorithm that computes each agent’s optimal allocation using only local computation and communication. Our algorithm resolves the temporal coupling through a local sweep, and the spatial, resource-sharing coupling through a distributed consensus iteration with every step feasibility guarantee. Under a mild diagonal-dominance assumption, we show that the optimal allocation at any fixed time converges geometrically to a stationary solution as the horizon grows, and that our distributed algorithm tracks this limit at a matching geometric rate. This guarantees a memory-efficient variant of the algorithm, in which each agent retains only a bounded window of past history, with the resulting error decaying exponentially in the window size. We validate these guarantees numerically, confirming the predicted convergence rates under dynamically changing team membership.

\subsection{Motivating Applications} 
We are motivated by applications of problem \eqref{problem:main} to continuous monitoring of unknown environments using multi-robot systems, such as wildfire monitoring or search and rescue of a lost person. In these applications, a team of robots must decide how to allocate its collective sensing or coverage effort across a large region. As a result, the robots divide the area into subareas, each surveyed during a given time slot, with the team’s attention shifting from one subarea to the next as the mission unfolds. This setting exhibits exactly the three structural features that motivate problem \eqref{problem:main}. First, team composition is not fixed: robots join and leave the active mission as they complete sub-tasks, recharge, or move in and out of communication range, so the number and identity of agents contributing to the allocation vary with time. Second, the shared resource constraint, the total sensing, communication, or coverage budget available to the team, must hold exactly at every instant. Third, each robot’s own allocation is temporally coupled. Because its current assignment affects its position, battery state, or task progress, its allocation at one time slot directly constrains what is achievable at the next. Thus, the current decision must be made with an eye toward its effect on future allocations.

\subsection{Related Works}
Problem~\eqref{problem:main} is related to the literature on streaming and time-varying optimization studied in \cite{hamam2022streaming}, which considers an unconstrained, centralized variant of problem~\eqref{problem:main}. This is distinct from the online convex optimization paradigm \cite{zinkevich2003}, where a single variable is updated at each step and no notion of smoothing the past arises, and from prediction-correction time-varying optimization \cite{simonetto2020survey}, where the optimization variable does not grow with time. In this paper, we extend the streaming least-squares framework of \cite{hamam2022streaming} to a \emph{networked, multi-agent} setting with (i) a team of robots whose active membership changes at every time step, and (ii) a shared equality constraint coupling all active robots at each instant.

A separate and complementary line of work studies \emph{distributed} resource allocation: dividing a fixed quantity of resource among a network of agents, each with a private convex cost, using only local computation and neighbor communication. The foundational approach here is the center-free, weighted-gradient method of Xiao and Boyd \cite{xiao2006optimal}, in which each node updates its allocation in proportion to the difference between its own marginal cost and those of its neighbors; starting from a feasible point, this update preserves the resource constraint exactly at every iteration, exactly as our inner consensus loop does in Section~\ref{sec:algorithm}. This line was later extended to {dynamic} (time-varying) communication graphs, most notably by Doan and Olshevsky \cite{doan2017distributed}, who give a deterministic distributed protocol whose convergence time scales only quadratically in the network size over an arbitrary sequence of time-varying graphs. These results, however, assume each node's cost is a \emph{static} function of its own allocation alone: the graph may vary with time, but the per-node decision variables and their costs do not evolve or couple across time steps.

Our work also relates to the growing literature on open multi-agent systems (OMAS), where agents may enter or leave the network over time. Existing results study OMAS stability through dimension-varying switched-system models and Lyapunov criteria \cite{Xue2022_M3D_OMAS}, as well as contractive or paracontractive operator frameworks for consensus-type dynamics \cite{Franceschelli2021_StabOMAS,Oliva2023_SumStates,Deplano2024_ParacontractiveOMAS}. Several consensus protocols have also been developed to handle agent arrivals and departures, including max-median, gossip, ratio-consensus, and acknowledgement-based schemes \cite{Abdelrahim2017_MaxConsensus,Hendrickx2017_GossipRandomArrivals,Dashti2022_ModeOMAS,Makridis2024_OpenRC}. None of these works, however, considers the temporal coupling across time steps that is central to our formulation, so their results do not carry over to the setting studied here. \vspace{0.1cm}

\noindent\textbf{Notation.} We denote by $\mathbf{1}$ a vector whose entries are $1$. Given a function $f(x_{t-1},x_{t})$, we denote by $\nabla_{1} f(\cdot,\cdot)$ and $\nabla_{2} f(\cdot,\cdot)$ its partial gradient w.r.t the first and second arguments, respectively. Finally, we present the classic forward-backward sweeping method to solve the equation $\mathbf{M}\textbf{x} = \textbf{g}$, where $\mathbf{M}$ is a tridiagonal matrix
\begin{align}
\mathbf{M} = \begin{bmatrix}
H_1 & E_1^T & & \\
E_1 & H_2 & \ddots & \\
 & \ddots & \ddots & E_{T-1}^T \\
 & & E_{T-1} & H_T,
\end{bmatrix},\label{notation:tridiagonal-matrix}
\end{align}
with $E_{i},H_{i}$ are block submatrix with proper dimensions. The LU factorization of $\mathbf{M}$ is $\mathbf{M} = \mathbf{V}\mathbf{U}$, where $\mathbf{V},\textbf{U}$ given as upper and lower diagonal matrix, respectively, 
\begin{align*}
\textbf{V} = \begin{bmatrix}
V_1  & & \\
E_1 & V_2  & \\
 & \ddots & \ddots \\
  & & E_{T-1} & V_T
\end{bmatrix},
\textbf{U} = \begin{bmatrix}
I & U_1  & & \\
 & \ddots & \ddots \\
 &  & I & U_{T-1} \\
 &  & & I
\end{bmatrix}.
\end{align*}
The pivots $V_t$ and multipliers $U_t$ are defined recursively as
\begin{align*}
V_1 = H_1, \quad V_t = H_t - E_{t-1}V_{t-1}^{-1}E_{t-1}^T,\quad
U_t = V_t^{-1}E_t^T.
\end{align*}
With the LU factorization, equation $\mathbf{M}\mathbf{x} = \mathbf{g}$ becomes $\mathbf{V}\mathbf{U}\mathbf{x} = \mathbf{g}$. Let $\mathbf{z}=\mathbf{U}\mathbf{x}$, $\mathbf{V}\mathbf{U}\mathbf{x} = \mathbf{g}$ splits into two triangular systems, $\mathbf{V}\mathbf{z}=\mathbf{g}$ followed by $\mathbf{U}\mathbf{x}=\mathbf{z}$. First, we consider $\mathbf{V}\mathbf{z}=\mathbf{g}$. Since $\mathbf{V}$ is upper diagonal matrix, we can apply a forward iteration to find $\mathbf{z}$ as follows \cite{golub2013matrix}
\begin{equation*}
z_1 = Q_1^{-1}g_1, \qquad z_t = V_t^{-1}\big(g_t - E_{t-1}z_{t-1}\big), \quad t=2,\dots,T.
\end{equation*}
Once we find $\mathbf{z}$, we apply a backward update, going from $T$ to $1$, to compute $\mathbf{x}$ as follows 
\begin{equation*}
x_T = z_T, \qquad x_t = z_t - U_t\,x_{t+1}, \quad t=T-1,\dots,1.
\end{equation*}

\section{Time-Varying Quadratic Optimization\\ in Open Networks}
We study problem \eqref{problem:main} where each function $f_{i,t}$ is quadratic
\begin{align}
\hspace{-0.3cm}f_{i,t}(x_{i,t-1},x_{i,t}) = \big(B_{i,t}x_{i,t-1}+A_{i,t}x_{i,t}-y_{i,t}\big)^2 + \gamma_{i} x_{i,t}^2,\label{objective:f_i}
\end{align}
where $y_{i,t},A_{i,t},B_{i,t}$ are constants and  $\gamma_{i}$ is a regularized parameter. Associating a scalar multiplier $\lambda_t$ with the resource constraint at each time $t$, the Lagrangian of \eqref{problem:main} is
\begin{equation}
\mathcal{L}(\mathbf{x}_{1:T},\boldsymbol{\lambda}_{1:T}) = \sum_{t=1}^{T}f_{t}(x_{t-1},x_{t}) + \lambda_t\Big(\mathbf{1}^Tx_{t} - b_t\Big).
\label{eq:lagrangian}
\end{equation}
Using $\mathcal{L}$, the optimality condition of \eqref{problem:main} is
\begin{align*}
    \left\{\begin{array}{l}
     \sum_{t=1}^{T}\nabla_{2}f_{t}(x_{t-1},x_{t}) + \lambda_{t}\mathbf{1} = 0,  \\
     \mathbf{1}^Tx_{t} - b_{t} = 0.
    \end{array}\right.
\end{align*}
For agent $i$ active at time $t$, its local cost depends on $x_{i,t}$ in two ways across the sum in \eqref{problem:main}: as the \emph{second} argument of $f_{i,t}(x_{i,t-1},x_{i,t})$, and as the \emph{first} argument of $f_{i,t+1}(x_{i,t},x_{i,t+1})$. Note that if $i$ is not active at time $t+1$ then we simply have $f_{i,t+1} = 0$. For convenience, we define
\begin{align*}
H_{i,t} &:= A_{i,t}^2 + B_{i,t+1}^2 + \gamma_i, & E_{i,t} &:= A_{i,t+1}B_{i,t+1},\\
g_{i,t} &:= A_{i,t}y_{i,t} + B_{i,t+1}y_{i,t+1}.
\end{align*}
Thus, using \eqref{objective:f_i}, we obtain
\begin{align*}
&\nabla_2 f_{i,t}(x_{i,t-1},x_{i,t}) + \nabla_1 f_{i,t+1}(x_{i,t},x_{i,t+1})\\
&= 2 H_{i,t}\,x_{i,t} + 2E_{i,t}\,x_{i,t+1} + 2E_{i,t-1}\,x_{i,t-1} - 2g_{i,t}, 
\end{align*}
which by using the optimality condition above gives for all $t$
\begin{equation*}
H_{i,t}\,x_{i,t} + E_{i,t}\,x_{i,t+1} + E_{i,t-1}\,x_{i,t-1} - g_{i,t} = -\tfrac{1}{2}\lambda_t.
\end{equation*}
This condition implies that at any time $t$,  agent $i$ achieves a consensus on their gradient w.r.t their variable $x_{i,t}$. We next present a matrix form to represent this condition. We introduce the following notation
\begin{equation*}
H_t := \mathrm{diag}\big(H_{i,t}\big), \;
E_t := \mathrm{diag}\big(E_{i,t}\big), \;
g_t := \big(g_{i,t}\big),\; {i\in\mathcal{V}_t},
\end{equation*}
which gives a matrix form of the consensus condition above 
\begin{equation*}
H_t\,x_t + E_t\,x_{t+1} + E_{t-1}\,x_{t-1}  - g_t =  -\tfrac{1}{2}\lambda_t\,\mathbf{1}.
\end{equation*}
Second, let $L_t$ denote the (weighted) graph Laplacian of $\mathcal{G}_t$. Since $\mathcal{G}_{t}$ is connected (Assumption \ref{assum:connectivity}), we have $L_t\mathbf{1}=0$, or $\ker L_t = \mathrm{span}\{\mathbf{1}\}$. Using $L_{t}$, we obtain 
\begin{equation*}
L_{t}(H_t\,x_t + E_t\,x_{t+1} + E_{t-1}\,x_{t-1}  - g_t) = \mathbf{0}.
\end{equation*}
Finally, define $\mathbf{g}_{1:T} = (g_{t})$, $\mathbf{L} = \mathrm{diag}(L_1,\dots,L_T)$, and the block-tridiagonal matrix $\mathbf{M}$ as in \eqref{notation:tridiagonal-matrix}. 
The optimality condition of problem \ref{problem:main} becomes 
\begin{equation}
\mathbf{L}\big(\mathbf{M}\mathbf{x}_{1:T} - \mathbf{g}_{1:T}\big) = \mathbf{0}\quad\text{and}\quad \mathbf{1}^Tx_{t} = b_{t},\; \forall t.
\label{optimality:Laplacian-consensus-overall}
\end{equation}

\section{Distributed Consensus Forward-Backward Sweeping Algorithms}
\label{sec:algorithm}
In this section, we present a distributed algorithm, based on a forward-backward sweeping update, to solve \eqref{optimality:Laplacian-consensus-overall}. Since $\mathbf{M}$ is block-tridiagonal, it admits an LU factorization $\mathbf{M}=\mathbf{VU}$, where $\mathbf{V}$ and $\mathbf{U}$ are, respectively, block lower- and upper-bidiagonal matrices. Substituting into \eqref{optimality:Laplacian-consensus-overall} gives
\begin{align*}
\mathbf{L}\big(\mathbf{VU}\mathbf{x}_{1:T} - \mathbf{g}_{1:T}\big) = \mathbf{0}.
\end{align*}
Introducing the intermediate variable $\mathbf{z}_{1:T}=\mathbf{U}\mathbf{x}_{1:T}$, we split the equation above into two triangular systems: a \emph{forward} system $\mathbf{L}(\mathbf{V}\mathbf{z}_{1:T} -\mathbf{g}_{1:T})=\mathbf{0}$, and a \emph{backward} system $\mathbf{U}\mathbf{x}_{1:T}=\mathbf{z}$. Our main idea is to solve the forward system with a distributed consensus iteration, and then recover $\mathbf{x}_{1:T}$ through a purely local backward update.

Because each agent's objective is temporally coupled only to itself, $\mathbf{V}$ and $\mathbf{U}$ are themselves block-diagonal across agents, and every agent can compute its own entries without any communication. At time $t$, agent $i\in\mathcal{V}_t$ locally computes
\begin{equation*}
V_{i,1}=H_{i,1}, \qquad V_{i,t} = H_{i,t} - \frac{E_{i,t-1}^2}{V_{i,t-1}}, \qquad U_{i,t} = \frac{E_{i,t}}{V_{i,t}}\cdot
\end{equation*}
The proposed distributed consensus forward-backward sweeping method is formally presented in Algorithm~\ref{alg:DCFBS}, whose two main components are distributed consensus forward updates and local backward sweeping. We discuss them in turn below.\vspace{0.1cm}

\noindent\textbf{Distributed consensus forward updates}, implemented in Steps 3--9 of Algorithm~\ref{alg:DCFBS}, mirror the classic distributed consensus-gradient method for solving a static networked resource allocation problem \cite{xiao2006optimal}. Here, $r_i^{(k)}$ is agent $i$'s local gradient at inner iteration $k$. Each agent exchanges its local gradient with its neighbors and forms a consensus update using the weights $L_{t,ij}$ of the Laplacian matrix $L_t$ (Step~8). The goal of this step is to approximately solve $L_t(V_tz_t-g_t)=0$ using $K$ distributed consensus-gradient steps, so that $z_{i,t}$ converges to a value on which every agent's local gradient value agrees.\vspace{0.1cm}

\noindent\textbf{Local backward sweeping.} After $K$ consensus updates, each agent $i$ sweeps backward in time to revise its earlier local variables $x_{i,\ell}$ (Step~12). This step uses only agent $i$'s own locally stored quantities $z_{i,\ell}$ and $U_{i,\ell}$, and requires no further communication between agents.

\begin{algorithm}[h]
\caption{Distributed Consensus Forward-Backward Sweeping Algorithm}
\label{alg:DCFBS}
\begin{algorithmic}[1]
\For{$t=1,\dots,T$}
    \ForAll{agents $i$}
    \State $\tilde g_{i,t}\gets g_{i,t} - E_{i,t-1}z_{i,t-1}$
    \State $z_i^{(0)}\gets b_t/|\mathcal V_t|$ for all $i\in\mathcal V_t$ 
    \For{$k=0,\dots,K-1$}
        \State $r_i^{(k)}\gets V_{i,t}z_i^{(k)}-\tilde g_{i,t}$
        \State Broadcast $r_i^{(k)}$ to neighbors in $\mathcal G_t$
        \State $z_i^{(k+1)}\gets z_i^{(k)} - \alpha\sum_{j\in\mathcal N_{i,t}}L_{t,ij}\big(r_i^{(k)}-r_j^{(k)}\big)$ 
    \EndFor
    \State $z_{i,t}\gets z_i^{(K)}$ 
    \State $x_{i,t}\gets z_i^{(K)}$
    \For{$\ell=t-1,\dots,1$}
        \State $x_{i,\ell}\gets z_{i,\ell} - U_{i,\ell}\,x_{i,\ell+1}.$\quad  \textit{(Backward sweeping)}
\EndFor
    \EndFor
\EndFor
\State \textbf{Output:} $\{x_{i,t}\}$
\end{algorithmic}
\end{algorithm}

\begin{remark}
\noindent\textbf{Memory requirement.} Because the backward sweep at time $T$ relies on $z_{i,\ell}$ and $U_{i,\ell}$ for every earlier $\ell=1,\dots,T-1$, each agent must retain its own local history of these quantities for the full horizon observed so far, giving a per-agent memory footprint that grows as $O(T)$. In the next section, we show that under a mild diagonal-dominance assumption, each agent instead needs to store only a bounded window of its past history, at the cost of a small, explicitly quantified loss in accuracy that decays exponentially in the window size.
\end{remark}

\section{Theoretical Results}
\label{sec:theory}
We now characterize the convergence behavior of Algorithm~\ref{alg:DCFBS} under the following assumptions. 

\begin{assumption}
For every $t$, $\mathcal{G}_t$ is connected and undirected.
\label{assum:connectivity}
\end{assumption}
\begin{assumption}
There exist $\kappa\ge1$, $\delta\in (0,1)$, $\theta\in (0,\frac{1-\delta}{2})$ such that, for all active $i,t$,
\begin{equation}
|\kappa^{-1}H_{i,t}-1|\le\delta, \qquad |E_{i,t}|\le\kappa\theta.
\label{eq:assum-loose}
\end{equation}
\label{assum:loose-revised}\vspace{-0.5cm}
\end{assumption}
This assumption basically represents a weakly temporal coupling condition of the objective function $f_{i}$ at different time step. To see this, let us recall that $H_{i,t} := A_{i,t}^2 + B_{i,t+1}^2 + \gamma_i,  E_{i,t} := A_{i,t+1}B_{i,t+1}.$ The condition $\theta\le(1-\delta)/2$ implies that the temporal coupling $E_{i,t}$ linking agent $i$'s allocation across consecutive time steps be small relative to its own diagonal curvature $H_{i,t}$, specifically $|E_{i,t}|\le H_{i,t}/2$. This is structurally the same requirement that underlies the classical convergence theory of Jacobi and Gauss--Seidel iterations for solving a linear system $Ax=b$ \cite{golub2013matrix}. In particular, these methods split $A=D-L-U$ into its diagonal and off-diagonal parts and iterate using only the diagonal's inverse, and both are guaranteed to converge, from any starting point, whenever $A$ is \emph{strictly diagonally dominant}, i.e.\ $|a_{ii}|>\sum_{j\ne i}|a_{ij}|$ for every row $i$. Diagonal dominance guarantees that the spectral radius of the resulting iteration matrix is strictly less than one, so that errors introduced anywhere in the system are damped, rather than amplified, as they propagate through successive iterations. The condition in \eqref{eq:assum-loose} is basically equivalent to this diagonal dominance condition.  

To study the convergence of Algorithm \ref{alg:DCFBS}, we denote by $x^{\star}_{i,t\mid T}$ the optimal solution of problem \eqref{problem:main}, where the time horizon $T$ is fixed. Our main result is presented as follows.

\begin{proposition}
Under Assumption~\ref{assum:loose-revised}, there exists a limit $x_{i,t}^*$ such that $\lim_{T\rightarrow\infty}x_{i,t\mid T}^{\star} = x_{i,t}^{\star}$, for all $i\in\mathcal{V}_{t}$ and $t\in [1,T]$. Moreover, $\{x_{i,t}\}$ generated by Algorithm \ref{alg:DCFBS} satisfies 
\begin{equation}
\| x_{i,t} - x_{i,t}^*\| \le C_{1}\,\rho^{\,T-t} + \frac{C_{2}}{1-\rho}\,\mu^K, \quad \rho = \frac{\theta}{1-\varepsilon_*} < 1,
\label{eq:geom-conv-revised}
\end{equation}
where constants $\varepsilon_*,C_{1},C_{2},$ are defined in \eqref{eq:eps-star-def} and \eqref{constant-def}, respectively. These constants depend on $\kappa,\theta,\delta,$ and the algebraic connectivity of the Laplacian $L_{t}$. 
\label{prop:convergence-revised}
\end{proposition}
First, Proposition~\ref{prop:convergence-revised} shows that, as the horizon $T$ grows, the time-varying optimization problem in \eqref{problem:main} admits a stationary solution at each time $t$, independent of $T$. This is not obvious a priori: each new time step appends a new cost term and a new set of variables to the problem, and in general there is no reason that revising the solution to account for this new information should ever stabilize, rather than continuing to drift indefinitely.

Second, Proposition~\ref{prop:convergence-revised} further shows that the solution produced by the distributed consensus forward-backward sweeping algorithm converges to this same stationary solution at a geometric rate as the horizon grows. The corresponding rate constant depends only on the problem's weak-coupling parameters, through $\kappa,\theta,\delta$, the network topology, and the number of inner steps $K$.

One potential drawback of Algorithm~\ref{alg:DCFBS} is that the backward sweep requires each agent to store its local quantities $z_{i,\ell}$ and $U_{i,\ell}$ for every time step up to $T$, which becomes memory-inefficient when $T$ is large. Proposition~\ref{prop:convergence-revised}, which shows that these updates converge geometrically, suggests that we can instead limit the number of backward steps without sacrificing too much of the accuracy. We formalize this by bounding the additional error incurred when the sweep is limited to a buffer of size $B$, which keeps both the memory footprint and the computational cost bounded by $B$ rather than $T$. The resulting truncated algorithm is a simple modification of Algorithm~\ref{alg:DCFBS}, in which the backward sweep in Steps~11--13 runs for only $B$ steps rather than back to $t=1$.
\begin{align*}
    \hspace{-0.8cm}\textbf{for } &\ell = 1,\ldots, B-1,\;\textbf{do} \\
   \hspace{-0.8cm} & \hat{x}_{i,t} = \gets z_{i,\ell} - U_{i,\ell}\,\hat{x}_{i,\ell+1}.\quad \textit{(Backward sweeping)}\\
    \hspace{-0.8cm}\textbf{end} &\textbf{ for}
\end{align*}
Let $\tilde{\mathbf{x}}_{t} = \hat{\mathbf{x}}_{t+1\mid t+B}$. The following corollary presents a bound on the difference between $\tilde{\mathbf{x}}_{t}$ and $\mathbf{x}^{\star}_{t}$. 
\begin{corollary}
Under the same conditions as in Proposition \ref{prop:convergence-revised}
\begin{equation}
\|\mathbf{x}^{\star}_{t} - \tilde{\mathbf{x}}_{t}\| \le C_{3}\,\gamma^{B},
\label{eq:buffer-bound}
\end{equation}
for all $0\leq t\leq T-B$, where $\gamma \in (0,1)$ and $C_{3}$ is some constant.
\label{cor:buffer-revised}
\end{corollary}
\begin{proof}
    The proof of this result follows the same argument as that of Proposition~\ref{prop:convergence-revised}, with the corresponding bounds adapted in the obvious way; we therefore omit it for brevity. 
\end{proof}

\section{Numerical Validation}
\label{sec:numerical}

In this section, we provide numerical simulations to validate the correctness of our algorithm and the theoretical results established in the previous sections. We organize the validation around three questions, addressed in turn: does $x_{i,t|T}$ actually converge to a stationary solution $x_{i,t}^\star$ as the horizon $T$ grows (Proposition~\ref{prop:convergence-revised})? Does it do so at the predicted geometric rate? And does truncating the backward sweep to a small, fixed buffer $B$ (Corollary~\ref{cor:buffer-revised}) incur only a small, controlled loss in accuracy?

\begin{figure}[h]
\centering
\includegraphics[width=1\columnwidth]{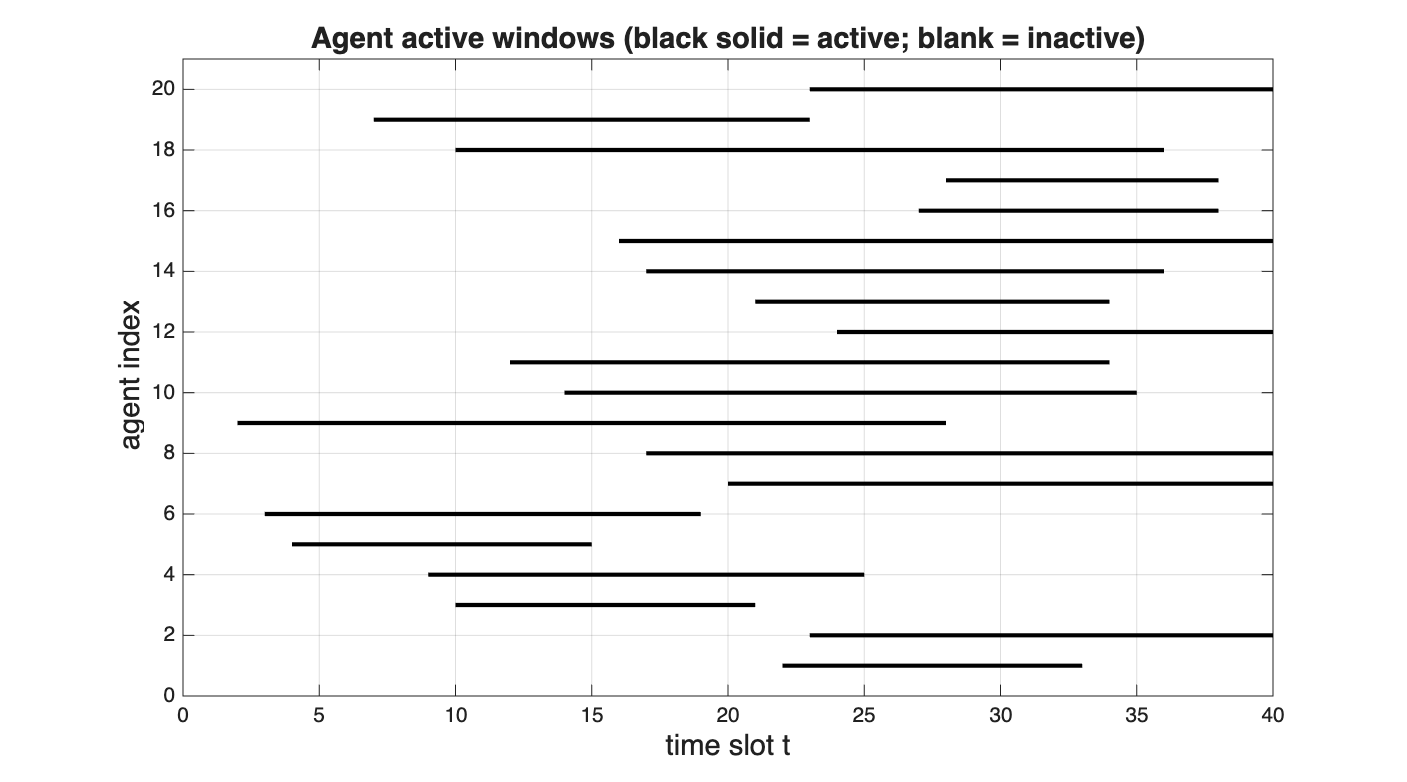}\vspace{-0.3cm}
\caption{Agent active profiles over horizon $T$.}
\label{fig:active}
\end{figure}

For our simulation, we consider a total of $20$ agents over a horizon of $T=40$ time slots, each with an active window within the horizon (Fig.~\ref{fig:active}). We enforce at least three agents active at every slot to keep each $\mathcal G_t$ non-trivial. Active agents communicate over a ring graph, with the Laplacian's off-diagonal entries between neighboring active agents set to $L_{t,ij}=-1/3$. Each agent's local matrices $A_{i,t},B_{i,t}$ are drawn as independent perturbations of $1$ and $0.6$, respectively, with regularization $\gamma_i\in[0.10,0.15]$. The coefficient $y_{i,t}$ is a smooth, agent-specific sinusoid plus small noise, and the shared budget $b_t$ scales with the active team size $n_t$ plus a slowly varying perturbation. We use $K=120$ consensus iterations per time step with step size $\alpha=0.35$. For an ease of exposition, we highlight three agents (agents $7$, $8$, and $15$), chosen randomly to show the results of our simulations.

\begin{figure}[h]
\centering
\includegraphics[width=1\columnwidth]{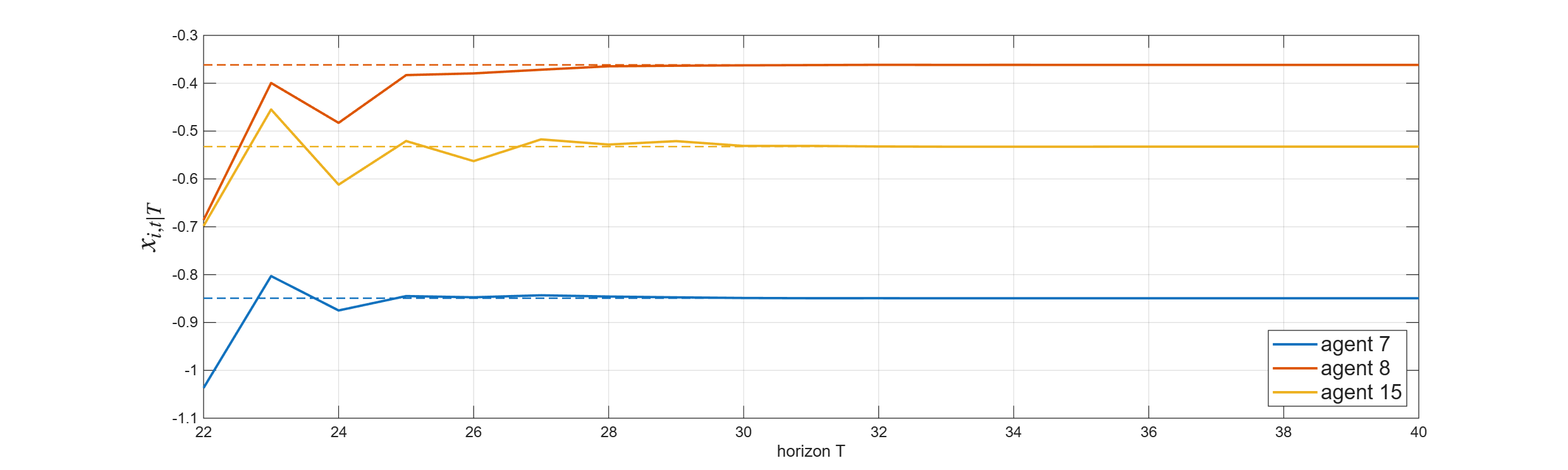}\vspace{-0.3cm}
\caption{Convergence of $x_{i,t|T}$ to $x_{i,t}^\star$ as $T$ grows.}
\label{fig:convergence_x_i_star}
\end{figure}

\noindent\textbf{Convergence of $x_{i,t|T}$ to $x_{i,t}^\star$} is shown in Fig.~\ref{fig:convergence_x_i_star}. We plot the value of the backward-swept estimate $x_{i,t|T}$, for agents $7$, $8$, and $15$, as the horizon $T$ is increased from $t$ up to $T_{\max}=40$, at three representative time slots ($t=14,20,26$) drawn from the agents' common active window. Dashed lines mark the corresponding values of $x_{i,t}^\star$, taken as $x_{i,t|T_{\max}}$. In every panel, each agent's trajectory visibly settles onto its dashed reference value well before $T$ reaches the horizon's end, with the largest fluctuations occurring immediately after $t$ and essentially disappearing within a number of additional time steps. These plots illustrate the first result of Proposition~\ref{prop:convergence-revised}, namely that $\lim_{T\rightarrow\infty}x_{i,t\mid T} = x_{i,t}^{\star}$.

\begin{figure}[h]
\hspace{-1cm}\includegraphics[width=1.2\columnwidth]{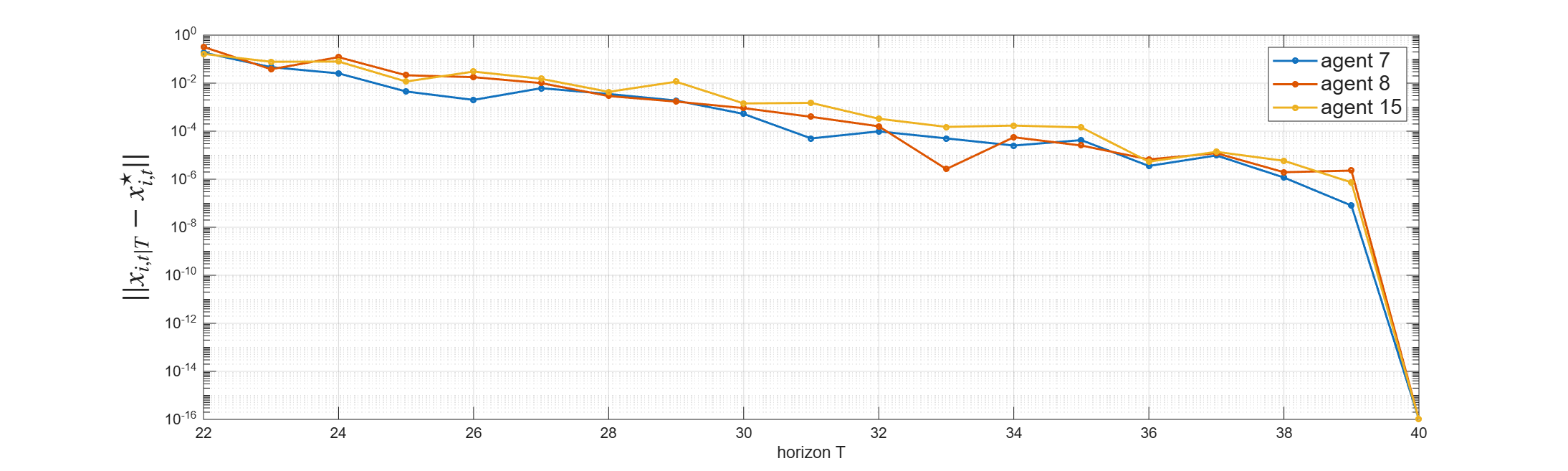}\vspace{-0.3cm}
\caption{Geometric convergence rate of $\|x_{i,t|T}-x_{i,t}^\star\|$ as $T$ grows.}
\label{fig:convergence_rate}
\end{figure}

\vspace{0.1cm}
\noindent\textbf{Geometric convergence rate of $\|x_{i,t|T}-x_{i,t}^{\star}\|$}. Fig.~\ref{fig:convergence_rate} shows the error $\|x_{i,t|T}-x_{i,t}^\star\|$ on a logarithmic scale, for the same agents and time slots as in Fig.~\ref{fig:convergence_x_i_star}. The goal is to show the quantitative rate in Proposition~\ref{prop:convergence-revised}. As the horizon $T$ increases, the error decays approximately linearly on this log scale, which agrees with the geometric decay with a constant rate $\rho<1$ in our theoretical results.

\begin{figure}[h]
\centering
\includegraphics[width=1\columnwidth]{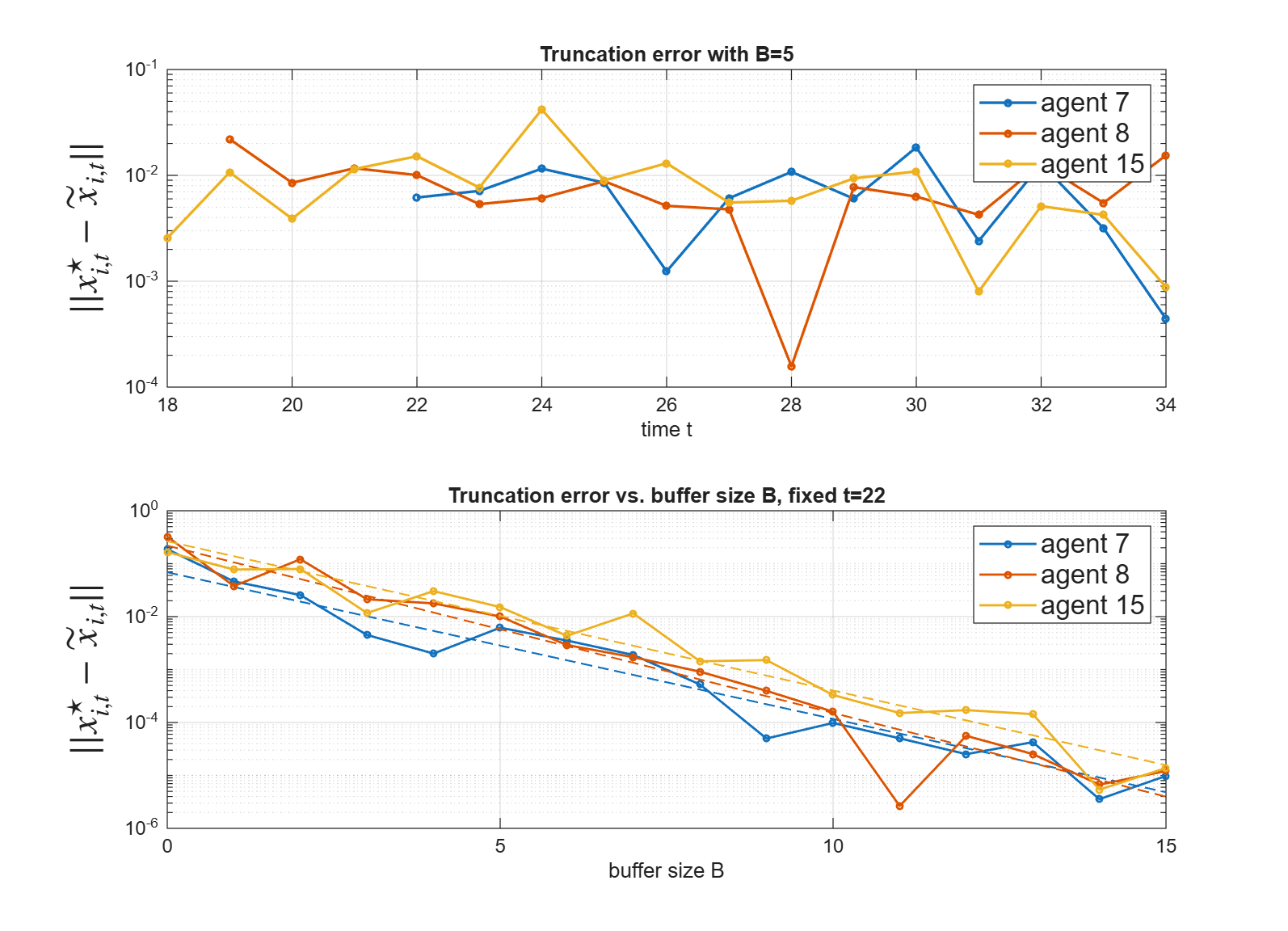}\vspace{-0.4cm}
\caption{Performance of Algorithm \ref{alg:DCFBS} under limited buffer $B$.}
\label{fig:finite-buffer}
\end{figure}

\noindent\textbf{Finite-buffer truncation $B$.} Finally, we investigate the approximation error introduced by finite-memory operation, as characterized in Corollary~\ref{cor:buffer-revised}. Specifically, we restrict the backward sweep to $B=5$ steps rather than the full $T=40$ steps used in Algorithm~\ref{alg:DCFBS}. Figure~\ref{fig:finite-buffer} (top) shows the resulting error $\|x_{i,t}^\star - \tilde x_{i,t}\|$ across different $t$, where $\tilde x_{i,t}$ denotes the truncated estimate returned with buffer size $B=5$. The error remains small and roughly uniform across $t$, confirming that limiting the backward sweep to just $5$ steps, rather than the full $40$, incurs little loss in accuracy. To verify the convergence rate predicted by Corollary~\ref{cor:buffer-revised}, Figure~\ref{fig:finite-buffer} (bottom) fixes the representative time $t=14$ and plots this error as a function of the buffer size $B$. As expected, the error decays to zero as $B$ increases, and does so linearly on a logarithmic scale, confirming the geometric convergence rate established in our theoretical result. Together, the two panels confirm that a small, horizon-independent buffer suffices for near-optimal accuracy.

\section{Conclusion}

We formulated streaming resource allocation over a multi-robot team with time-varying membership as a sequence of locally-coupled quadratic programs linked by a per-instant shared resource constraint, showed that the resulting bordered block-tridiagonal KKT system reduces to a scalar marginal-cost consensus condition when robots have no pairwise coupling, and derived a fully distributed algorithm combining consensus-gradient spatial updates with a temporal forward--backward smoothing sweep. We established geometric convergence and a finite-memory error bound, and validated the approach in simulation with dynamically changing team membership. Future work includes extending the framework to nonlinear robot dynamics via a distributed Newton-type update, relaxing the reliance on a connected communication graph at every instant, and validating the approach on a physical multi-robot testbed under realistic communication constraints.


\bibliographystyle{IEEEtran}
\bibliography{references}

\appendix
In this section, we present the proofs of our main results in Proposition \ref{prop:convergence-revised} and Corollary \ref{cor:buffer-revised}. We will assume that Assumptions \ref{assum:connectivity} and \ref{assum:loose-revised} always hold. Let $\mu_{t}\in(0,1)$ be the second smallest eigenvalue of the Laplacian $L_{t}$. 

Recall that the main idea of Algorithm \ref{alg:DCFBS} is to solve two systems of equations: a \emph{forward} system $\mathbf{L}(\mathbf{V}\mathbf{z}_{1:T} -\mathbf{g}_{1:T})=\mathbf{0}$, and a \emph{backward} system $\mathbf{U}\mathbf{x}_{1:T} = \mathbf{z}_{1:T}.$ Indeed, the inner loop of Algorithm is to find $\mathbf{z}_{1:T}$, while the backward sweeping update is to compute $\mathbf{x}_{1:T}$.

We first state the convergence of the inner loop for the consensus algorithm, which is studied in \cite{doan2017distributed}. Let $z_{i,t}$ be the solution of the forward equation and $z_{i}^{(K)}$ be the output of the consensus step after $K$ iterations. Then we have
\begin{align}
    |z_{i,t} - z_{i}^{(K)}| \leq  C_{t}\mu_{t}^{K}\label{result:consensus},
\end{align}
where $C_{t}$ is a constant depending on $1-\mu_t$ and the coefficients of the local objective functions. Let $C:=\sup_tC_t < \infty$ and $\mu:=\sup_t\mu_t<1$.

\subsection{Proof of Proposition~\ref{prop:convergence-revised}}\label{sec:proof-prop1}
We will leverage the techniques in  \cite{doan2017distributed} and \cite{hamam2022streaming} to derive the results in Proposition \ref{prop:convergence-revised}. 
Recall that 
\begin{align*}
H_{i,t} &:= A_{i,t}^2 + B_{i,t+1}^2 + \gamma_i,\quad  E_{i,t} := A_{i,t+1}B_{i,t+1},\\
g_{i,t} &:= A_{i,t}y_{i,t} + B_{i,t+1}y_{i,t+1},\\
V_{i,1}&=H_{i,1}, \qquad V_{i,t} = H_{i,t} - \frac{E_{i,t-1}^2}{V_{i,t-1}}, \qquad U_{i,t} = \frac{E_{i,t}}{V_{i,t}},
\end{align*}
where $V_{i,t}$ and $U_{i,t}$ are the LU factorization of the tridiagonal matrix $\mathbf{M}$. Let $M_y:=\sup_t|g_{i,t}|<\infty$ for all $i,t$. We observe that for any fixed $t$, the quantities $V_t,U_t,v_t$ do \emph{not} depend on the horizon $T$, provided $T\ge t+1$. 

To show the results in Proposition \ref{prop:convergence-revised}, we will first established a sequence of intermediate results as follows. First, we show that $V_{i,t}$ is bounded for all $i,t$.

\begin{lemma}
Under Assumption~\ref{assum:loose-revised}, $|\kappa^{-1}V_{i,t}-1|\le\varepsilon_*$ for all $i$ and $t\ge1$, where
\begin{equation}
\varepsilon_* = \frac{1+\delta}{2} - \sqrt{\left(\frac{1-\delta}{2}\right)^2-\theta^2} \;<\;1.
\label{eq:eps-star-def}
\end{equation}
\label{lem:pivot-bound}
\end{lemma}

\begin{proof}
Let $q_{i,t}:=\kappa^{-1}V_{i,t}$. We prove, by induction that 
\begin{equation*}
|q_{i,t}-1|\le\varepsilon_t \qquad\text{and}\qquad \varepsilon_t<1 \quad \forall t,
\end{equation*}
where $\varepsilon_1:=\delta$ and $\varepsilon_t:=\delta+\theta^2/(1-\varepsilon_{t-1})$ for $t\ge2$. 

\emph{Base case} ($t=1$): Since $V_{i,1}=H_{i,1}$, by Assumption~\ref{assum:loose-revised} and $\varepsilon_1=\delta<1$ we have
$$|q_{i,1}-1|=|\kappa^{-1}H_{i,1}-1|\le\delta=\varepsilon_1,$$
which gives $q_{i,1}\ge1-\delta>0$.

\emph{Inductive step:} Suppose $|q_{i,t-1}-1|\le\varepsilon_{t-1}<1$. This implies
\begin{equation*}
q_{i,t-1} \;\ge\; 1-\varepsilon_{t-1} \;>\; 0.
\end{equation*}
Dividing the recursion $V_{i,t}=H_{i,t}-E_{i,t-1}^2/V_{i,t-1}$ by $\kappa$,
\begin{equation*}
q_{i,t} - 1 = \big(\kappa^{-1}H_{i,t}-1\big) - \frac{(E_{i,t-1}/\kappa)^2}{q_{i,t-1}}.
\end{equation*}
By Assumption~\ref{assum:loose-revised}, we have $|\kappa^{-1}H_t-1|\le\delta$ and $|E_{t-1}/\kappa|\le\theta$. Thus, using $q_{t-1}\ge1-\varepsilon_{t-1}>0$ gives
\begin{equation*}
|q_{i,t}-1| \;\le\; \delta+\frac{\theta^2}{q_{i,t-1}} \;\le\; \delta+\frac{\theta^2}{1-\varepsilon_{t-1}} \;=\; \varepsilon_t.
\end{equation*}
It remains to show $\varepsilon_t<1$. Our approach is to show $\varepsilon_t\le\varepsilon_* < 1$ for all $t$. By Assumption \ref{assum:loose-revised}, $\theta\leq \frac{1-\delta}{2}$, which gives 
\begin{equation*}
\varepsilon_* = \frac{1+\delta}{2} - \sqrt{\left(\frac{1-\delta}{2}\right)^2-\theta^2} \;<\;1.
\end{equation*}
Note that $\varepsilon_{\star}$ satisfies
\begin{align*}
\varepsilon_*=\delta+\frac{\theta^2}{1-\varepsilon_*} \geq \delta. 
\end{align*}
We next show that $\epsilon_{t}\leq \epsilon_{\star}$ by induction. 

\emph{Base case} ($t=1$): We have $\varepsilon_*\ge\delta=\varepsilon_1$.

\emph{Inductive step:} suppose $\varepsilon_{t-1}\le\varepsilon_* < 1$. Consider
\begin{equation*}
\varepsilon_* - \varepsilon_t = \frac{\theta^2}{1-\varepsilon_*} - \frac{\theta^2}{1-\varepsilon_{t-1}}= \frac{\theta^2}{(1-\varepsilon_*)(1-\varepsilon_{t-1})}\big(\varepsilon_*-\varepsilon_{t-1}\big).
\end{equation*}
Since $\varepsilon_{t-1} \leq \varepsilon_{\star} < 1$, the right-hand side is nonnegative. Thus,  $\varepsilon_t\le\varepsilon_*<1$, which completes our proof. 
\end{proof}


\begin{lemma}
We have $|U_{i,t}|\le\rho:=\theta/(1-\varepsilon_*)<1$ for all $t$.
\label{lem:U-contraction}
\end{lemma}

\begin{proof}
By Assumption \ref{assum:loose-revised} and Lemma~\ref{lem:pivot-bound}, we have 
$$|U_{i,t}|=\frac{|E_{i,t}|}{|V_{i,t}|}\le\frac{\kappa\theta}{\kappa(1-\varepsilon_*)}=\rho.$$ 
By the definition of $\varepsilon_{\star}$ we have
$$\varepsilon_*=\delta+\frac{\theta^2}{1-\varepsilon_*},$$ which gives $\theta^2=(\varepsilon_*-\delta)(1-\varepsilon_*)$. Thus, by \eqref{eq:eps-star-def} we have 
\begin{equation*}
\rho^2 = \frac{\theta^2}{(1-\varepsilon_*)^2} = \frac{\varepsilon_*-\delta}{1-\varepsilon_*} < 1,
\end{equation*}
which concludes our proof.
\end{proof}

\begin{lemma}
We have $|z_{i,t}|\le M_z$ and $|x_{i,t}|\le M_x$ for all $t\le T$, where
\begin{equation}
M_z := \frac{M_y}{\kappa(1-\varepsilon_*)(1-\rho)}, \qquad M_x := \frac{M_z}{1-\rho}.
\label{eq:Mv-Mx-def}
\end{equation}
\label{lem:uniform-bound}
\end{lemma}

\begin{proof}
Recall from the forward update that 
\begin{align*}
    z_{i,t} = \frac{g_{i,t} - E_{i,t-1}z_{i,t-1}}{V_{i,t}},
\end{align*}
which by using  Lemma~\ref{lem:pivot-bound} gives
\begin{equation*}
|z_{i,t}| \le \frac{|g_{i,t}|+|E_{i,t-1}||z_{i,t-1}|}{V_{i,t}} \le \frac{M_y}{\kappa(1-\varepsilon_*)} + \rho\,|z_{i,t-1}|.
\end{equation*}
Using a standard induction argument, we obtain
\begin{align*}
    |z_{i,t}| \le \frac{M_{y}}{\kappa(1-\varepsilon_{\star})(1-\rho)}\cdot  
\end{align*}
To show the bound of $x_{i,t}$ we will apply an inductive proof.  For the base case, $| x_{i,T}|=|z_{i,T}|\le M_z\le M_x$ holds. Inductively, for $t<T$, using Lemma~\ref{lem:U-contraction},
\begin{equation*}
|x_{i,t}| \le |z_{i,t}| + |U_{i,t}||x_{i,t+1}| \le M_z + \rho M_x = \frac{M_z}{1-\rho},
\end{equation*}
which concludes our proof.
\end{proof}

\begin{proof}[Proof of Proposition~\ref{prop:convergence-revised}]
Let $x_{i,t\mid T}$ be the solution returned by Algorithm \ref{alg:DCFBS} under horizon $T$, and recall $x_{i,T\mid T}=z_i^{(K)}$ and $x_{i,t\mid T}=z_i^{(K)}-U_{i,t}x_{i,t+1\mid T}$ for $t<T$, where $z_i^{(K)}$ is computed once, when time $t$ is first processed, and does not depend on $T$. Using this recursion, we consider for all $t=T-1,T-2,\dots$
\begin{equation*}
x_{i,t|T+1} - x_{i,t|T} = -U_{i,t}\big(x_{i,t+1|T+1}- x_{i,t+1|T}\big),
\end{equation*}
which by using $|U_{i,t}|\le\rho$ (Lemma~\ref{lem:U-contraction}), we obtain
\begin{equation*}
\big|x_{i,t|T+1} - x_{i,t|T}\big| \le \rho^{\,T-t}\,\big|x_{i,T|T+1}- x_{i,T|T}\big|.
\end{equation*}
At $t=T$, the preceding relation does not apply directly since $x_{i,T|T}=z_i^{(K)}$. By the triangle inequality and Lemma~\ref{lem:uniform-bound},
\begin{equation*}
\big|x_{i,T|T+1}- x_{i,T|T}\big| \le \big| x_{i,T|T+1}\big|+\big| x_{i,T|T}\big| \le 2M_x,
\end{equation*}
which gives
\begin{equation*}
\big|x_{i,t|T+1}- x_{i,t|T}\big| \le 2M_x\,\rho^{\,T-t}, \qquad \forall\,t\le T.
\end{equation*}
To show the convergence of $x_{i,t\mid T}$, we will show it is a Cauchy sequence. For $N>T\ge t$, we consider
\begin{equation*}
x_{i,t|N} - x_{i,t|T} = \sum_{S=T}^{N-1}\big( x_{i,t|S+1}- x_{i,t|S}\big),
\end{equation*}
which gives
\begin{align}
&\big|x_{i,t|N} - x_{i,t|T} \big| \le 2M_{x}\sum_{S=T}^{N-1}\rho^{\,S-t}\notag\\
&\le 2M_x\rho^{\,T-t}\sum_{\ell=0}^{\infty}\rho^\ell = \frac{2M_x}{1-\rho}\,\rho^{\,T-t} \xrightarrow[T\to\infty]{} 0,\notag
\end{align}
uniformly in $N>T$. Hence $\{x_{i,t|T}\}_{T\ge t}$ is Cauchy in $T$, therefore it is a convergent sequence. Let the limit point of $\{x_{i,t|T}\}$ be $x_{i,t}^{\star(K)}$. Then we have
\begin{equation}
\big|x_{i,t|T}-x_{i,t}^{\star(K)}\big| \le \frac{2M_x}{1-\rho}\,\rho^{\,T-t}.
\label{eq:step1-fixed-K}
\end{equation}
We next relate $x_{i,t}^{\star(K)}$ to the true limit $x_{i,t}^\star$ of Proposition~\ref{prop:convergence-revised}. Let $z_{i,t}$ denote the exact solution of the forward equation at time $t$, and let $\tilde{x}_{i,t\mid T}$ denote the \emph{idealized} backward sweep obtained by using $z_{i,t}$ in place of $z_i^{(K)}$ at every step, i.e. $\tilde{x}_{i,t\mid T}=z_{i,T}$ and $\tilde{x}_{i,t\mid T}=z_{i,t}-U_{i,t}y_{i,t+1\mid T}$ for $t<T$. Since $z_{i,t}$ likewise does not depend on $T$, the same argument as above applies verbatim to $\{\tilde{x}_{i,t\mid T}\}_{T\ge t}$, which is therefore also Cauchy in $T$ and converges to $x_{i,t}^\star$, with
\begin{equation}
\big|x_{i,t|T}-x_{i,t}^\star\big| \le \frac{2M_x}{1-\rho}\,\rho^{\,T-t}.
\label{eq:step1-idealized}
\end{equation}
Using $x_{i,t\mid T}=(z_i^{(K)}-z_{i,t})+z_{i,t}-U_{i,t}x_{i,t+1\mid T}$ together with $\tilde{x}_{i,t\mid T}=z_{i,t}-U_{i,t}\tilde{x}_{i,t+1\mid T}$, we consider $e_{i,t}:=x_{i,t\mid T}-\tilde{x}_{i,t\mid T}$, which satisfies $e_{i,T}=z_i^{(K)}-z_{i,T}$ and, for $t<T$,
\begin{equation*}
e_{i,t} = \big(z_i^{(K)}-z_{i,t}\big) - U_{i,t}\,e_{i,t+1}.
\end{equation*}
Using $|z_{i,t}-z_i^{(K)}|\le C_t\mu_t^K\le C\mu^K$ (Lemma~\ref{lem:U-contraction}), we obtain
\begin{equation*}
|e_{i,t}| \le C\mu^K + \rho\,|e_{i,t+1}|,
\end{equation*}
and telescoping from $t$ up to $T$ exactly as in Lemma~\ref{lem:uniform-bound} gives, uniformly in $T\ge t$,
\begin{equation*}
\big|x_{i,t\mid T}-\tilde{x}_{i,t\mid T}\big| = |e_t| \le C\mu^K\sum_{\ell=0}^{T-t}\rho^\ell \le \frac{C}{1-\rho}\,\mu^K.
\end{equation*}
Letting $T\to\infty$, and using $x_{i,t|T}\to x_{i,t}^{\star(K)}$ together with $\tilde{x}_{i,t|T}\to x_{i,t}^\star$, we obtain
\begin{equation*}
\big|x_{i,t}^{\star(K)}-x_{i,t}^\star\big| \le \frac{C}{1-\rho}\,\mu^K,
\end{equation*}
which when combining with \eqref{eq:step1-fixed-K} gives
\begin{align*}
\big| x_{i,t|T}-x_{i,t}^*\big| &\le \big|x_{i,t|T}-x_{i,t}^{\star(K)}\big| + \big|x_{i,t}^{\star(K)}-x_{i,t}^\star\big|\notag\\ 
&\le \frac{2M_y}{\kappa(1-\varepsilon_*)(1-\rho)^3}\rho^{\,T-t} + \frac{C}{1-\rho}\,\mu^K.
\end{align*}
This concludes our proof where 
\begin{align}
    C_{1} := \frac{2M_y}{\kappa(1-\varepsilon_*)(1-\rho)^3},\quad C_{2} := \frac{C}{1-\rho}\cdot\label{constant-def} 
\end{align}
\end{proof}

\end{document}